\documentclass[a4paper]{paper}
\usepackage{fullpage}
\usepackage{graphicx,amsfonts,latexsym,color,hyperref}

\usepackage{algorithmic}
\usepackage{algorithm}
\usepackage{amsmath}

\usepackage{epsfig}
\usepackage{makeidx}
\usepackage{amssymb}
\usepackage{latexsym}
\usepackage{amsfonts}

\newcommand{\Real}[1]{ { {\mathbb R}^{#1} } }

\newcommand{\be}{\begin{equation}}
\newcommand{\ee}{\end{equation}}
\newcommand{\bea}{\begin{eqnarray}}
\newcommand{\eea}{\end{eqnarray}}
\newcommand{\bean}{\begin{eqnarray*}}
\newcommand{\eean}{\end{eqnarray*}}
\newcommand{\ba}{\begin{array}}
\newcommand{\ea}{\end{array}}
\newcommand{\bit}{\begin{itemize}}
\newcommand{\eit}{\end{itemize}}
\newcommand{\ben}{\begin{enumerate}}
\newcommand{\een}{\end{enumerate}}

\newtheorem{lemma}{Lemma}

\newtheorem{proposition}{Proposition}

\newcommand{\beai}{\begin{eqnarray}}
\newcommand{\eeai}[1]{\label{eq:#1}\end{eqnarray}}
\newcommand{\nn}{\nonumber \\}
\newcommand{\ese}{\end{eqnarray*}}
\newcommand{\bse}{\begin{eqnarray*}}
\newcommand{\rf}[1]{~(\ref{eq:#1})}
\def\Argmin{\mathop{\hbox{\rm Argmin$\,$}}}

\begin{document}
\title{Robust Eigenvector of a Stochastic Matrix with Application to PageRank}
\author{Anatoli Juditsky  and Boris Polyak%
\thanks{A. Juditsky is with LJK, Universit\'e J. Fourier, BP 51, 38041 Grenoble Cedex 09, France {\tt juditsky@imag.fr}}
\thanks{B. Polyak is with Institute for Control Sciences, 65 Profsoyuznaya str., 117997 Moscow, Russia {\tt
boris@ipu.ru}}}
\date{\today}
\maketitle \thispagestyle{empty}
\begin{abstract}
We discuss a definition of robust dominant eigenvector of a family
of stochastic matrices.
Our focus is on application to ranking problems, where the proposed approach can be seen as a robust
alternative to the standard PageRank technique.
The robust eigenvector computation is reduced to a convex
optimization problem. We also propose a simple algorithm for robust eigenvector approximation which can be viewed
as a regularized power method with a special stopping rule.
\end{abstract}

\section{Introduction}
\label{S:intro}
There are numerous problems which can be reformulated as finding a
dominant eigenvector of a stochastic matrix --- computing invariant
probabilities in Markov chains, various problems in econometrics,
sociometry, bibliometrics, sport etc., see an excellent survey on
history of such applications \cite{giants}. Recently, ranking in
the Web (PageRank) attracted a lot of attention \cite{BrinPage,LanMeyer}.

Let us briefly remind this approach and discuss some problems which arise in relation to it. There
are $n$ pages (nodes of a graph) connected with outgoing and
incoming links. Then $ x_i$ --- score of i-th page of the Web,
$i=1,...,n$ is defined as follows:
 {\em $x_i = \sum_{j\in L_i} x_j/n_j$, where $n_j$ is the number of outgoing links of page $j$ and $L_i$
 is the set of pages linked to page $i$}.
 Thus the score vector $x$ satisfies the equation
 \beai
 Px=x
 \eeai{1up}
  with $P_{ij}=1/n_j$ if page $j$ cites page $i$
 and $P_{ij}=0 $ otherwise. I.e. $x$ is an eigenvector
 corresponding to eigenvalue $1$ of the (column) stochastic matrix $P$ of size $n\times n$.

In the example below $n=7$, outgoing links are denoted by arrows.
 \begin{figure}[htb]
\centerline{\includegraphics[width=0.5\linewidth]{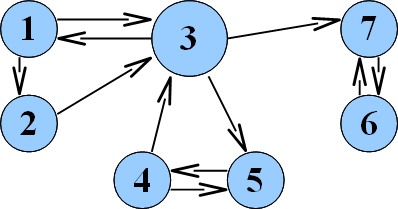}}
\caption{Test example}\label{test}
\end{figure}
The corresponding matrix $P$ and its dominant (or principal) eigenvector $\bar{x}$ are

 $$P=\left(\begin{tabular}{ccccccc}
  0 & 0 & 1/3 & 0 & 0 & 0 & 0 \\
  1/2 & 0 & 0 & 0 & 0 & 0 & 0 \\
  1/2 & 1 & 0 & 1/2 & 0 & 0 & 0 \\
  0 & 0 & 0 & 0 & 1 & 0 & 0 \\
  0 & 0 & 1/3 & 1/2 & 0 & 0 & 0 \\
  0 & 0 & 0 & 0 & 0 & 0 & 1 \\
  0 & 0 & 1/3 & 0 & 0 & 1 & 0 \\
\end{tabular}\right),$$
$$\bar{x}=(0,0,0,0,0,0.5,0.5)^T.$$

We see that the scores of pages 6 and 7 equal 0.5 while all other
pages have zero scores. Of course it contradicts common sense and our
intuition about ``true scores''. There are several problems related to the above
definition:
\begin{enumerate}
\item We have assumed that $n_j\neq 0$ for all $j$. Yet, in applications, there are nodes
  with no outgoing links --- \emph{dangling nodes} (e.g. pdf files).
  One can easily tackle this problem by taking $n_j=n$ for all dangling
  nodes. However more complicated problem of ``traps'' (absorbing nodes 6 and 7 in the above example) remains open.
\item In the case of disconnected
  subgraphs (typical for real-life problems) $\bar{x}$ is not unique.
\item  Small changes in $P$ may result in dramatic perturbations of $\bar{x}$ ---
  the computed scores are non-robust.
\item  Power method with iterations $x_{k+1}=Px_k$ do not converge for cyclic
  matrices; for acyclic ones it may converge slowly due to the small gap $1-|\lambda_2|$
  (here  $\lambda_2$ is the second eigenvalue of $P$).
\end{enumerate}
PageRank medicine for these problems goes back to the pioneering
paper \cite{BrinPage}: allow the $i$-the page to receive a part of its score ``for free'' -- replace
the dominant eigenvector $\bar{x}$ of $P$ with the (unique) solution to the equation
\[x^{\alpha}=\alpha Px^{\alpha}+(1-\alpha)e,\; e_i=1/n, \;i=1,...,n.\]
In other words, the ``true'' adjacency matrix $P$ is replaced with
\beai
M=\alpha P+(1-\alpha)E, \;\; E_{ij}=1/n, \;i,j=1,...,n,
\eeai{PR}
so that $x^{\alpha}=Mx^{\alpha}$ is the unique dominant eigenvector of stochastic matrix $M$.
Unlike a dominant eigenvector of $P$, $x^{\alpha}$ is
robust to variations in links, with fast convergence of power iterations (as
$\alpha^k$). As a justification the model of \emph{surfer in the Web}
is proposed in \cite{BrinPage}: a surfer follows links with probability $\alpha$ and jumps to
an arbitrary page with probability $1-\alpha$. Then $x^{\alpha}_i$
is the (asymptotically) average time share he spends on page $i$.
\par
Here we propose a different modification of the original score
definition \rf{1up}. Note that in most of applications nominal
matrices are known approximately and they are subject to
perturbations, so robustness issues become important. Our goal is to
provide a robust counterpart to standard techniques in the spirit of
Robust Optimization framework  \cite{BTN}.  In the hindsight, it
shares many common features with with Robust Least Squares
\cite{EGL}. To the best of our  knowledge this approach to ranking
problems is new.

\section{Problem formulation}
We say that matrix $P\in \Real{n\times n}$ is \emph{stochastic}
(column-stochastic) if $P_{ij}\geq 0 \quad\forall i,j, \quad\sum_i
P_{ij}=1 \quad\forall j$. The set of all stochastic matrices is
denoted $\mathcal{S}$.
The celebrated Perron-Frobenius theorem states that
there exist a \emph{dominant} eigenvector $\bar{x}\in \Sigma$ (here $\Sigma=\big\{u\in \Real{n}|\,\sum_i u_i=1,\, u_i\ge 0\big\}$
being the standard simplex of $\Real{n}$) of the ``nominal'' matrix $P$:
$$P\bar{x}=\bar{x}.$$
\par
We are looking for a robust extension of the dominant eigenvector.
Let $\mathcal{P}\subset  \mathcal{S}$ stand for the set of
\emph{perturbed} stochastic matrices $F=P+\xi$, where $\xi$ is a
perturbation. The function $\max_{F\in \mathcal{P}} \|Fx-x\|$, where
$\|\cdot\|$ is some norm, can be seen as  a measure of ``goodness''
of a vector $x$ as common dominant eigenvector of the family
$\mathcal{P}$. We say that the vector $\hat{x}$ is a robust solution
of the eigenvector problem on $\mathcal{P}$ if \beai
    \hat{x}\in \Argmin_{x\in \Sigma}\left\{\max_{F\in \mathcal{P}} \|Fx-x\|\right\}.
\eeai{def}
Let us consider some choices of the parameters of the above definition -- uncertainty set $\mathcal{P}$ and the norm $\|\cdot\|$ -- for the PageRank problem.
Recall that our motivation is to ``immunize'' the score $\hat{x}$  against small perturbations of the adjacency matrix $P$. We may consider,
that the uncertainty of the $j$-th column of $P$, in terms of the $\ell_1$-norm of elements, is bounded by $\varepsilon_j$.
In other words, if $[\xi]_j$ is the $j$-th column of the perturbation matrix, then $\|[\xi]_j\|_1\le \varepsilon_j$, $j=1,...,n$. From now on we assume that $\varepsilon_j>0$, $j=1,...,n$. Further, we may fix the ``total uncertainty budget'' $\varepsilon\ge \|\xi\|_1$, where with some notational abuse, $\|\xi\|_1=\sum_{ij}|\xi_{ij}|$. With a ``natural choice'' of the norm $\|\cdot\|=\|\cdot\|_1$, the definition \rf{def} now reads
\beai
    \hat{x}^{(1)}&\in &\Argmin_{x\in \Sigma}\Big\{\max_{\xi} \|(P+\xi)x-x\|_1\Big| \nn
    &&e^T[\xi]_{j}=0,\;\|[\xi]_j\|_1\le \varepsilon_j,\;\|\xi\|_1\le \varepsilon\Big\}
\eeai{all1}
(here $e_i=n^{-1}$, $i=1,...,n$).
Let $n_j$ be the number of outgoing links on page $j$. If we set $\varepsilon_j\le {1\over n_j}$, then $F=P+\xi$ is a stochastic matrix.
\par
A different robust eigenvector may be obtained if, instead of fixing a bound for the $\ell_1$-norm of the perturbation, we bound its Frobenius norm: $\|
\xi\|_F\le \varepsilon$. When combining it with the choice of the norm $\|\cdot\|=\|\cdot\|_2$ in \rf{def}, we come to the definition
\beai
   \hat{x}^{(2)}&\in &\Argmin_{x\in \Sigma}\Big\{\max_{\xi} \|(P+\xi)x-x\|_2\Big| \nn
    &&e^T[\xi]_{j}=0,\;\|[\xi]_j\|_1\le \varepsilon_j,\;\|\xi\|_F\le \varepsilon\Big\}.
\eeai{1and2}
An interesting simple version of the problem \rf{1and2} may be obtained when substituting the ``column-wise constraints'' on $\xi$ with the requirement that $P+\xi$ is stochastic:
\beai
   \hat{x}^{(F)}&\in&\Argmin_{x\in \Sigma}\Big\{\max_{\xi} \|(P+\xi)x-x\|_2\Big|\nn
   && F+\xi \mbox{ is stochastic, }\|\xi\|_F\le \varepsilon\}.
\eeai{all2}

Obviously, there is a vast choice of uncertainty sets and
norms, which will result in different  definitions of robust eigenvector.
Three definitions above may be viewed as a starting point for this study.
\section{Theoretical analysis}
Observe that the min-max problems involved in the definitions  \rf{all1} -- \rf{all2} of the robust eigenvector are not convex-concave. Indeed, if $\Xi$ is the corresponding set of perturbation matrices, the function
\[
\phi_\Xi(x)= \max_{\xi\in \Xi} \|(P+\xi)x-x\|
\]
 cannot be computed efficiently.
  However, one can easily  upper bound the optimal values of these problems by those of simple convex problems as follows.
\begin{proposition}
\label{spto1}
Let \[\Xi_1=\{\xi \in \Real{n\times n}\big|\;e^T[\xi]_{j}=0,\;\|[\xi]_j\|_1\le \varepsilon_j,\;\|\xi\|_1\le \varepsilon\},\] and let $\|\cdot\|=\|\cdot\|_1$. Then
\[
\phi_{\Xi_1}(x)\le \varphi_1(x),\;\varphi_1(x)=\|Px-x\|_1+{\varepsilon}g_1(x),
\]
where the convex function $g_1$ is defined as
\bse
g_1(x)=\min_{x=u+v}\Big\{ \|u\|_\infty+\sum_j{\varepsilon_j\over \varepsilon}|v|_j\Big\}.
\ese
\end{proposition}
The proofs of the statements of this section are postponed till the appendix.
\par
Note that $g_1$ is a norm. Further, when $\varepsilon_1=...=\varepsilon_n$ and $s={\varepsilon\over \varepsilon_1}$ is integer, $g_1(x)$ is simply the
sum of $s$ largest in magnitude elements of $x$. It is worth to mention that computing $g_1(x)$ requires $O(n)$ (up to a logarithmic factor in $n$) arithmetic operations.

We have an analogous bound for function $\phi_\Xi(x)$, involved in the definition \rf{1and2}:
\begin{proposition}
\label{spto2}
Let \[\Xi_2=\{\xi\in  \Real{n\times n}\big|\;e^T[\xi]_{j}=0,\;\|[\xi]_j\|_1\le \varepsilon_j,\;\|\xi\|_F\le \varepsilon\},\] and let $\|\cdot\|=\|\cdot\|_2$. Then
\[
\phi_{\Xi_2}(x)\le \varphi_2(x),\;\varphi_2(x)=\|Px-x\|_2+{\varepsilon}g_2(x),
\]
where the convex function $g_2$ satisfies
\bse
g_2(x)=\min_{x=u+v}\Big\{ \|u\|_2+\sum_j{\varepsilon_j\over \varepsilon}|v|_j\Big\}.
\ese
\end{proposition}
What was said about the numerical cost of computing $g_1(x)$ may be repeated here: computing $g_2(x)$ requires $O(n)$ (up to ``logarithmic factors'') arithmetic operations.
\par
We have a particularly simple and intuitive bound for \rf{all2}:
\begin{proposition}
\label{spto22}
Let  \[
\Xi_F=\{P+\xi \mbox{ is stochastic, }\|\xi\|_F\le \varepsilon\}.\] Then
\[\phi_{\Xi_F}(x)\leq \varphi_F(x), \; \varphi_F(x)=\|Px-x\|_2+\varepsilon \|x\|_2.\]
\end{proposition}
 Note that the upper bound $\varphi_F(x)$ for $\phi_{\Xi_F}$
is not very conservative. Indeed, the corresponding worst-case perturbation is
$$\xi^*=\frac{\varepsilon (Px-x)x^T}{\|Px-x\|\|x\|},$$
and for $F=P+\xi^*$ we have $\sum_i F_{ij}=1$, however the
condition $F_{ij}\geq 0$ (required for $F$ to be stochastic) may be
violated.
\par
The results above can be summarized as follows: \\
let $\|\cdot\|_{(1)}$ and $\|\cdot\|_{(2)}$ be some norms (in the above examples we have $\|\cdot\|_{(1)}$ set to $\ell_1$- or $\ell_2$-norm, and
$\|\cdot\|_{(2)}=g_1,\;g_2$ or $\|\cdot\|_2$). With certain abuse of terminology, we refer to the vector
\beai
    \hat{x}\in \Argmin_{x\in \Sigma}\left\{\varphi(x)=\|Px-x\|_{(1)}+\varepsilon \|x\|_{(2)}\right\}
\eeai{robeig}
 as (computable) \emph{robust dominant eigenvector} of the corresponding family of matrices.

Let us discuss some properties of the robust eigenvector.
\begin{enumerate}
  \item $\hat{x}$ coincides with $\bar{x}$ if $P$ is regular and
  $\varepsilon$ is small enough. Recall that stochastic
matrix $P$ is called \emph{regular}, if its largest eigenvalue
$\lambda_1=1$ is simple, while all other eigenvalues lay inside the
unit circle $1-|\lambda_i|\geq \mu>0, \;i=2,...,n$. For instance,
matrices with all positive entries are regular. Thus dominant
eigenvector of a regular matrix is robust with respect to small
perturbations. This is the case of matrix $M$ in \rf{PR}.
   \item In the case $\|x\|_{(2)}=\|\cdot\|_2$ vector $\hat{x}$
   is unique due to the strict convexity of $\|x\|_2$ on $\Sigma$.
\item For large $\varepsilon$ vector $\hat{x}$ is close to $e\in \Sigma, \;e_i=n^{-1},\; i=1,...,n$.
\end{enumerate}

\section{Solving optimization problem \rf{robeig}}\label{solving}
For irregular matrices or for matrices with small gap $\mu>0$ the robust
dominant vector  $\hat{x}$ can not be found in explicit form. To
compute it one should minimize the function $\varphi(x)$ over
standard simplex $\Sigma$. In the above examples these are well structured convex optimization problems, and as such they can be efficiently solved using available optimization software. For instance,  medium-size problems
(say, when $n\leq 1.e3-1.e4$) may be solved using interior-point methods. For large-scale problems one can use methods of mirror-descent family \cite{NY, LNS, NaPo09CDC} and solve a
saddle-point reformulation of \rf{robeig}, as well as various randomized algorithms (see
 \cite{IshTemp, NaPo09CDC} and references therein). Finally, for
huge-dimensional problems Yu. Nesterov \cite{Nesterov} proposed recently special randomized
algorithms which can work with sparse matrices and $n=1.e6-1.e9$. We do not discuss these methods here.
Instead, we propose a ``fast and dirty'' method for minimization of \rf{robeig},
which happens to be a close relative of the PageRank technique.
\par
The method starts the
minimization process from the point $x_1$, which is the minimizer
of $\|x\|_{(2)}$ over $\Sigma$. Then we apply the power method with averaging to minimize
the first term in $\varphi(x)$:
\[x_k=\frac{x_1+Px_1+...+P^{k-1}x_1}{k}\]
(note that power method without averaging may be non convergent,
for instance for a cyclic $P$).  We have
\[
\|Px_k-x_k\|_{(1)}=\frac{\|P^k x_1-x_1\|_{(1)}}{k}\leq \frac{c}{k},
\]
that is $\|Px_k-x_k\|_{(1)}=O(1/k)$.
\par
To fix the ideas, let us consider the case of Proposition \ref{spto22},
where the norms $\|\cdot\|_{(1)}$ and $\|\cdot\|_{(2)}$ are the Euclidean norms. We have $x_1=e$, and
when writing the method in recurrent form
we get
\begin{equation}\label{alg}
    x_{k+1}=\left(1-\frac{1}{k+1}\right)Px_k+\frac{1}{k+1}e.
\end{equation}
When denoting $M_k=\left(1-\frac{1}{k+1}\right)P+\frac{1}{k+1}E, \quad E_{ij}=1/n$,
the iterative method reads $x_{k+1}=M_k x_k, \quad x_1=e.$ We
apply the method until $\varphi_k=\varphi(x_k)$ starts to arise.
Thus we arrive to
\\{\em Algorithm 1}
\\{\tt begin:} $x_1=e$
\\{\tt $k$-th iteration:} $x_{k+1}=M_k x_k$,
\[
M_k=(1-\frac{1}{k+1})P+\frac{1}{k+1}E\]
\\{\tt stop:} $\varphi_{k+1}>\varphi_k, \bar{x}=x_k$\\
and $\bar{x}$ is the approximate solution. Note that the numerical complexity
of one iteration is dominated with that of multiplication $Px$ (and, of
course, we should not write $M_k$ in dense form for sparse $P$).

If we compare Algorithm 1 with PageRank algorithm: \[ x_{k+1}=M x_k, \quad
M=\alpha P+(1-\alpha) E, \quad\alpha=const,
\]
we conclude that our
algorithm is PageRank with varying $\alpha$, equipped with special
stopping rule. The proposed method is also closely related to the the regularized algorithm of \cite{ifac}. However, unlike the method of  \cite{ifac}, Algorithm 1 the stopping rule is related to the uncertainty level.
\par\noindent
{\bf Choice of $\varepsilon$} is an important issue for the described
approach. 
One may consider the following heuristics for Web ranking problems:
we may assume that the number
 $n_j$ of outgoing links of the $j$-th page is known with accuracy $\pm 1$ for $qn$ pages, for some $q<1$. Keeping in mind that the average $n_j=m=20$
for the whole Web, and taking $q=0.5$ we get for uncertainty in
$P$, measured in Frobenius norm, $\varepsilon\simeq
\sqrt{qn}/m\simeq 0.03\sqrt{n}=30$ for $n=10^6$. In small problems
(like the examples in this paper) the choice $\varepsilon\simeq 1$ appears to be convenient.

\section{Numerical results}
\label{numsec}
\subsection{Simulated graphs}
\begin{figure}[htb]
\centerline{\includegraphics[width=0.75\linewidth]{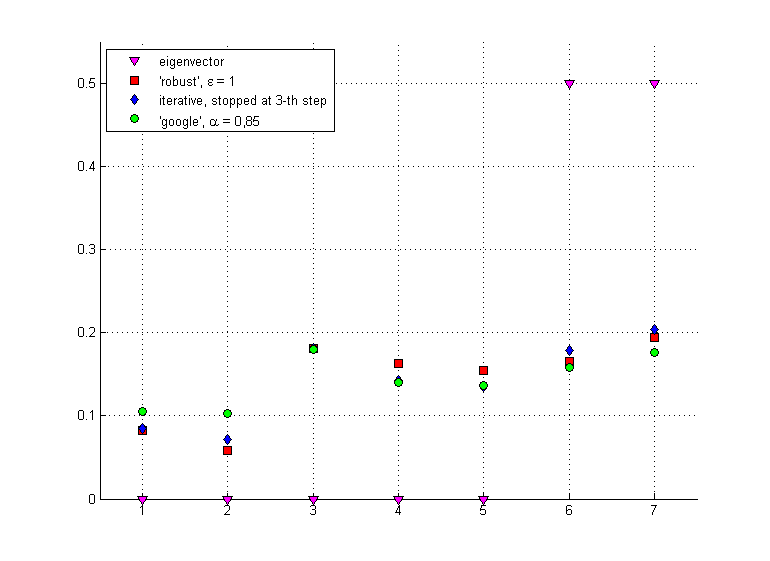}}
\caption{Results for test example} \label{results}
\end{figure}

To illustrate an application of the robust eigenvector, proposed above, let us go back to the example in the introduction. We present of Figure \ref{results} the scores computed for this example using  three algorithms: 1) PageRank algorithm with $\alpha=0.85$, 2) exact minimization of $\varphi_F$ with $\varepsilon=1$, 3) Algorithm 1 with $\varepsilon=1$.

We observe that the scores for all three methods are close and
differ dramatically from the dominant eigenvector of the
nominal matrix. The robust approach generates the scores which are acceptable from the ``common sense'' point of view. It is worth
to mentioning that Algorithm 1 requires only 4 iterations
($\varphi_{4}>\varphi_3, \bar{x}=x_3$).
\par
 We applied the proposed approach to  the
simulated examples proposed in \cite{ifac}. The corresponding graphs may be constructed in arbitrary dimensions (from modest to huge
ones), with the dominant eigenvectors of matrices
$P$ and $M$ which can be computed explicitly. Furthermore, the corresponding matrix-vector product $Px$ can be computed very efficiently without storing the
matrix $P$ itself.
\par
\begin{figure}[h]
\begin{center}
\begin{minipage}[h]{0.5\linewidth}
\includegraphics[width=1\linewidth]{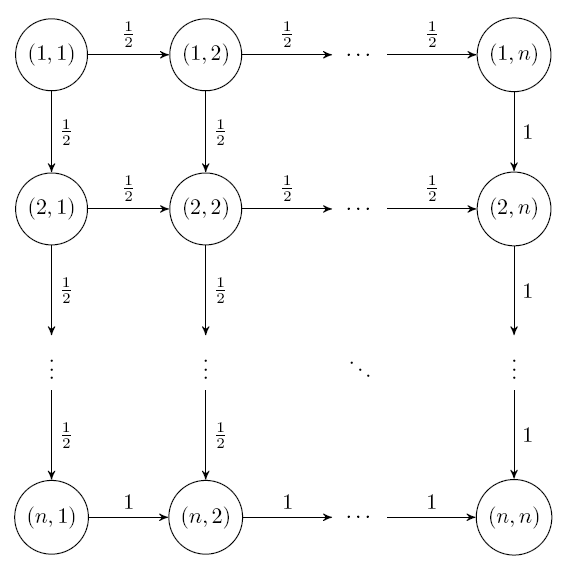}
\caption{Model 1}\label{Model 1}
\end{minipage}
\end{center}
\end{figure}

Consider the graph \emph{(Model 1)} of Figure \ref{Model 1}.
The nodes are shown in circles, the references are denoted by arrows. To get
rid of \((n,n)\) as a dangling node we assume the ability of equally
possible jump from this node to any other. There are $N=n^{2}$ nodes
in the model. We are interested in finding the vector $x=(x_{ij})$ of scores.

Taking arbitrary value of $x_{1 1}$, say $x_{1 1}=1$, we get the
system of equations equivalent to \rf{1up}:
$$\left\{
\begin{array}{ll}
    x_{i 1}=x_{1 i}=\frac{1}{2}x_{i-1, 1}+1,\quad i=2, ..., n,\quad x_{11}=1,\nonumber \\
    x_{i j}=\frac{1}{2}x_{i-1, j}+\frac{1}{2}x_{i, j-1}+1,\quad j, i=2, ..., n-1,\nonumber \\
    x_{n j}=x_{j n}=x_{n, j-1}+\frac{1}{2}x_{n-1, j}+1,\quad j=2, ..., n-1,\nonumber \\
    x_{n n}=x_{n-1, n}+x_{n, n-1}+1.\nonumber \\
\end{array}
\right.
$$
These equations can be solved explicitly (from $x_{1 1}$ we find
$x_{2 1}$ etc). Then we can proceed to standard normalization and
get the scores \(x_{ij}^*=\frac{x_{ij}}{\sum_{i,j}x_{ij}}\).

Similarly the scores for PageRank matrix $M$ can be found from
equations (it is convenient to take $y_{1 1}=\alpha$, then \mbox{$y_{n
n}=N\alpha$}):
$$\left\{
\begin{array}{ll}
    y_{i 1}=y_{1 i}=\alpha(\frac{1}{2}y_{i-1, 1}+1),\quad i=2, ..., n,\nonumber\\
    y_{i j}=\alpha(\frac{1}{2}y_{i-1, j}+\frac{1}{2}y_{i, j-1}+1),\qquad j, i=2, ..., n-1,\nonumber \\
    y_{n j}=y_{j n}=\alpha (y_{n, j-1}+\frac{1}{2}y_{n-1, j}+1),\qquad j=2, ..., n-1,\nonumber \\
    y_{n n}=\alpha (y_{n-1, n}+ y_{n, n-1}+1)\nonumber \\
\end{array}
\right.
$$
and normalization \(y_{ij}^*=\frac{y_{ij}}{\sum_{i,j}y_{ij}}\).
\par
Note that even very small perturbations may strongly
corrupt the graph on Figure \ref{Model 1}. Thus, taking large penalty parameter $\varepsilon$, as it was suggested in the discussion at the end of Section \ref{solving} would lead to
complete ``equalization'' of the scores. That is why we choose small
penalty coefficient in this example (specifically, it was taken
$\varepsilon=0.01$).
\par
\begin{figure}[htb]
\begin{center}
\includegraphics[width=0.75\linewidth]{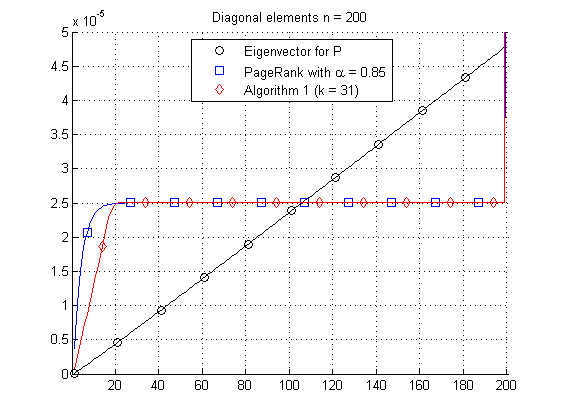}
\caption{Diagonal scores} \label{model1-diag}
\end{center}
\end{figure}
\begin{figure}[htb]
\begin{center}
\includegraphics[width=0.75\linewidth]{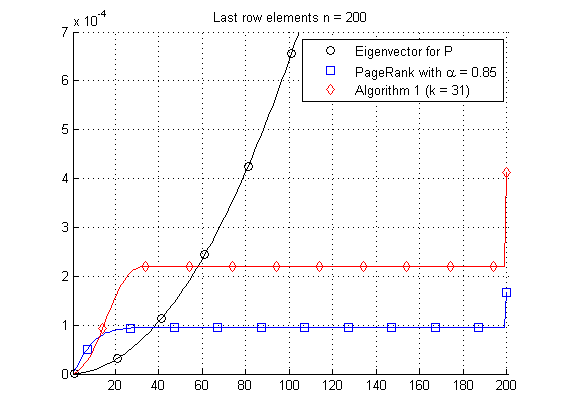}
\caption{Last row scores} \label{model1-lastrow}
\end{center}
\end{figure}

Let us denote $x_{\alpha}=y^*$. We compare now the score vector $x^*$ for
nominal $P$ with the score $x_{\alpha}$ by PageRank and the
results obtained by Algorithm 1. On Figure \ref{model1-diag} we
present the scores of diagonal nodes, i.e. $x_{i, i}, i=1,\dots ,
n$ for $n=200$ (that is $N=40 000$) and
$\alpha=0.85$.
Figure \ref{model1-lastrow} shows the scores $x_{n i}, i=1,\dots , n$ of the last row.
From the geometry of our graph (Figure \ref{Model 1}) it is natural to expect
the ``true'' scores $x_{ij}$ to increase with $i, j$ and the
highest scores should be achieved in south-east corner of the square.
The score vector, obtained as dominant eigenvector of $P$, possesses these
properties. PageRank solution and robust solution provide both
strongly ``equalized'' ranks; the robust solution appears to better fit the expected one.
\par
Let us modify  slightly the graph of Model 1: in the new \emph{Model
2} we link the node \((n, n)\) {\em only} to the node \((1, 1)\).
For this model the system of equations for which defines the nominal
scores is similar to the system for  Model 1, and the dominant
eigenvector of $P$ and PageRank vector $x^{\alpha}$ can be computed
explicitly. On the other hand, it is obvious that the power method
would not converge for this model. Indeed if $x^0_{1 1}=1, x^0_{i
j}=0, \{i j\}\neq \{1 1\}$ then after $2n-1$ iterations we return at
the same state, and matrix $P$ is cyclic.

\begin{figure}[htb]
\begin{center}
\begin{minipage}[h]{\linewidth}
\includegraphics[width=0.75\linewidth]{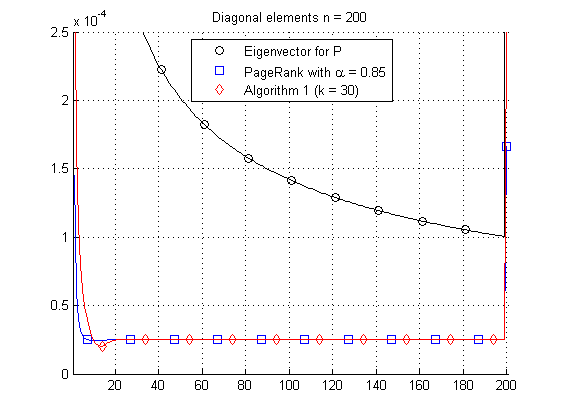}
\caption{Diagonal scores for the second model.} \label{model2-diag}
\end{minipage}
\hfill
\begin{minipage}[h]{\linewidth}
\includegraphics[width=0.75\linewidth]{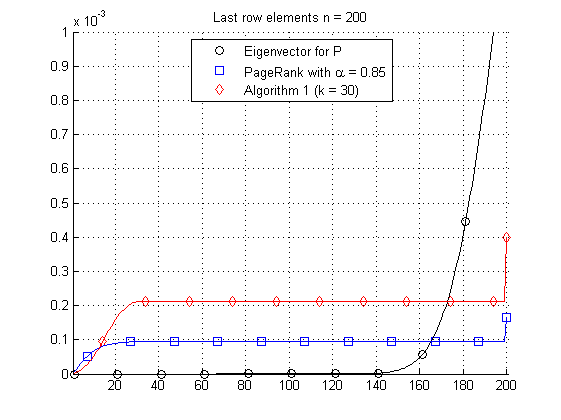}
\caption{Last row  scores for the second model.}
\label{model2-lastrow}
\end{minipage}
\end{center}
\end{figure}

The results for Model 2 with the same parameters ($n=200,
\varepsilon =0.01$) are presented at Figures \ref{model2-diag} and \ref{model2-lastrow}.
In this example, same as for Model 1, the score vector obtained
using PageRank is close to the robust one. And may be robust scores
are slightly closer to our expectations of ``common sense'' ranks.

\subsection{``Real-life'' data}
\begin{figure}[htb]
\begin{center}
\includegraphics[width=0.75\linewidth]{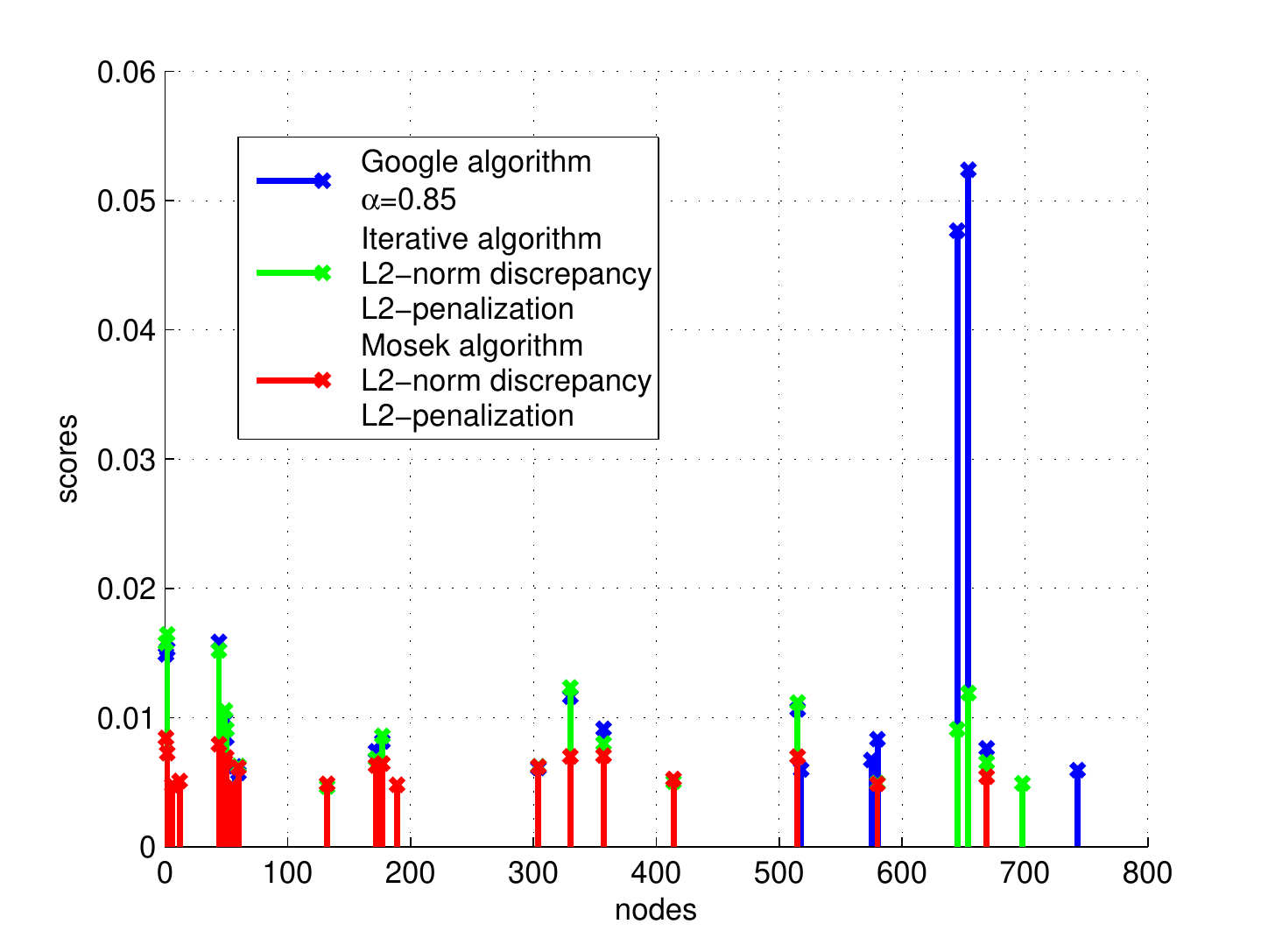}
\caption{20 highest scores for {\tt CompCom} matrix} \label{compcom}
\end{center}
\end{figure}

We have tested the proposed techniques on the experimental datasets, used in experiments on link analysis ranking algorithms described in \cite{www10}. We present here experimental results for adjacency matrices from two refined datasets:
{\tt Computational Complexity (CompCom)} and {\tt Movies}, downloaded from \href{http://www.cs.toronto.edu/\~tsap/experiments/}{http://www.cs.toronto.edu/\~tsap/experiments/}. {\tt CompCom}  data represent 789 pages with totalize 1449 links (with average number of links per page is equal to 1.8365). {\tt Movies} dataset has 4754 notes and 17198 (average number of links per page is  3.6176). We compare high accuracy solution to the problem \rf{robeig} (with the norms $\|\cdot\|_{(1)}$ and $\|\cdot\|_{(2)}$ are Euclidean norms) obtained using {\em MOSEK} optimization software \cite{Mosek} with the scores  by Algorithm 1 and those obtained by the PageRank algorithm. The same value of penalization parameter $\varepsilon=1$ in \rf{robeig} is used in all experiments. On Figure \ref{compcom} we present the 20 largest elements of the score vectors for the {\tt CompCom} data. In this experiment the optimal value of the problem \rf{robeig} is $0.0587$. The corresponding value for solution supplied by Algorithm 1 (after $3$ iterations) is $0.0756$.
It should be noted that the high scores $x^{\alpha}_{645}$ and $ x^{\alpha}_{654}$ by PageRank are pure artifacts -- they correspond to the pages which receive numerous links but only reference each other.
The techniques proposed in this paper seem better handle this sort of problem than the classical PageRank algorithm.
\par
\begin{figure}[htb]
\begin{center}
\includegraphics[width=0.75\linewidth]{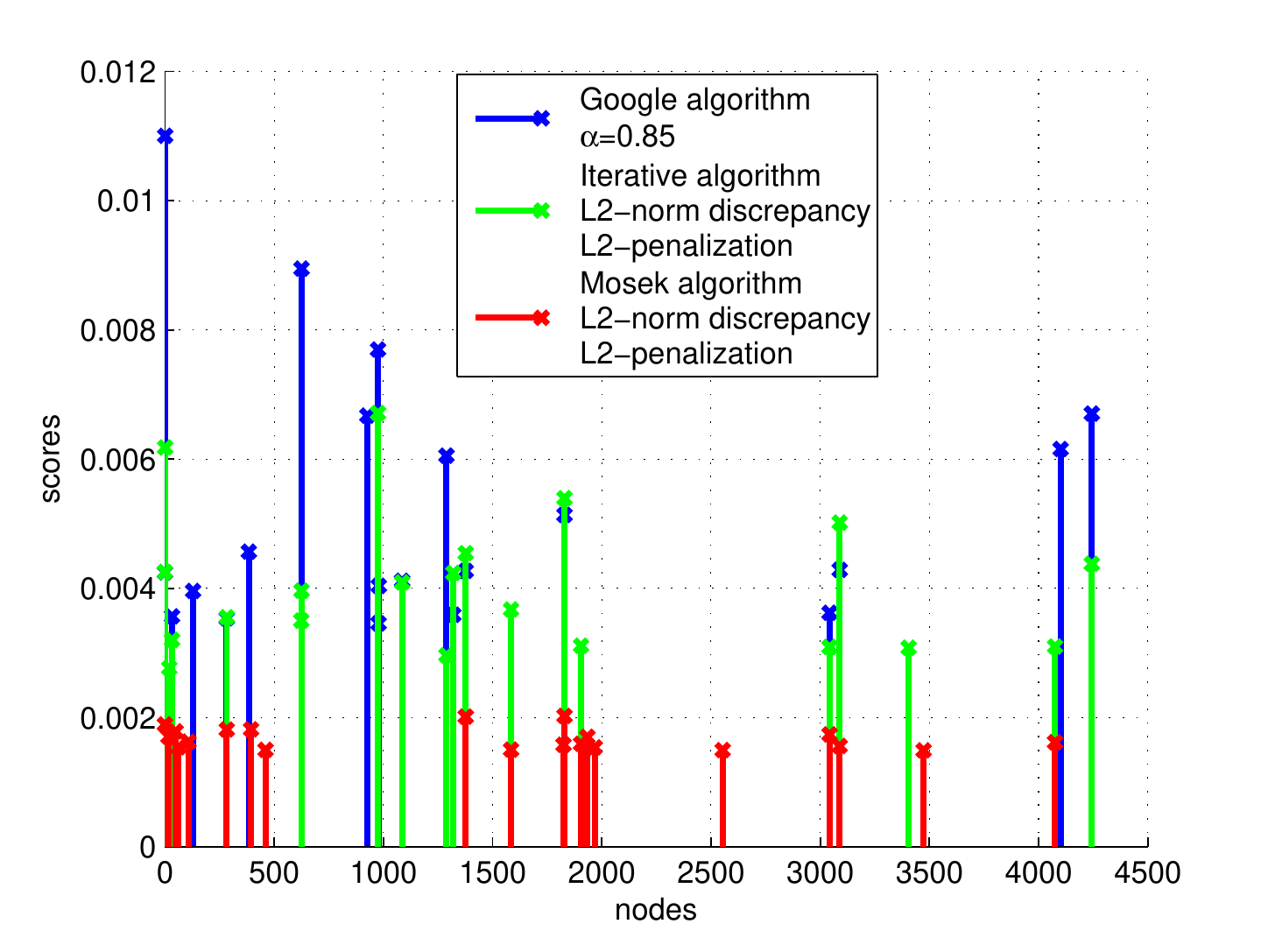}
\caption{20 highest scores for {\tt Movies} matrix} \label{movies}
\end{center}
\end{figure}

Figure \ref{movies} shows the 20 largest elements of the score vectors for the {\tt Movies} data. In this case the best scores computed by 3 competitors are quite close. The optimal value of the problem \rf{robeig} is $0.0288$; the approximate solution supplied by Algorithm 1 (after $4$ iterations) attains the value $0.0379$.

\section{Directions for future research}

\begin{enumerate}
  \item \emph{Structured perturbations.} It is often  natural
  to assume that the perturbation $\xi$  of nominal matrix $P$ possesses some structure. For instance, we
  can consider $\xi=M\Delta N$, where  $\|\Delta\|_F\leq
  \varepsilon$, and $M, N$ are rectangular matrices of
  appropriate dimensions, or $\xi=\sum t_i A_i$, where $\|t\|_2\leq
  \varepsilon$, with fixed matrices $A_i\in \Real{n\times n}$.
  \item \emph{Implicit uncertainty model.}
Family $\mathcal{P}$ can be defined implicitly. Suppose that we can sample randomly from $\mathcal{P}$, so that \emph{random} realizations
$F_i\in \mathcal{P},\; i=1,...,N$ are available. One may try to find $x$ which is approximate dominant
eigenvector for all of them. For instance, one can fix
$\varepsilon>0$ and look for solution of a system of convex inequalities
\[\|F_i x-x\|\leq \varepsilon, \quad i=1,...N, \quad x\in \Sigma.\]
Such system can be solved iteratively using the mirror-descent
algorithm. This approach is close to \cite{ineq, PolTempo}.
  \item \emph{Computational research.} The development of
   numerical methods tailored for optimization problem \rf{robeig}
  and their comparison is of great interest.
  \item \emph{Real-life examples.} We plan to test extensively the proposed approach on the large-scale Web rank problems data available from different
  Internet databases.
\end{enumerate}

\section{ACKNOWLEDGMENTS}
The authors would like to acknowledge Andrey Tremba's assistance when implementing numerical simulations of Section \ref{numsec} and enlightening discussions with Arkadi Nemirovski and Yuri Nesterov.
\appendix

\subsection{Proof of Propositions \ref{spto1} and \ref{spto2}}
Note that
\bse
\phi_{\Xi_1}(x)&=&\max_{\xi\in \Xi_1}\|(\xi+P-I)x\|_1\\
&=&\max_{\xi\in \Xi_1}\max_{u\in \Real{n},\|u\|_\infty\le 1}u^T(\xi+P-I)x.
\ese
Then $z=\xi^Tu$ satisfies
\bse
|z_j|&=&\left|\sum_{i} \xi_{ij}u_i\right|\le \sum_{i} |\xi_{ij}|\le \varepsilon_j.\\
\|z\|_1&=&\sum_{j}\left| \sum_{i} \xi_{ij}u_i\right|\le  \sum_{j}\sum_{i}|\xi_{ij}|\le
\varepsilon.
\ese
We conclude that
\beai
\phi_{\Xi_1}(x)&\le& \varphi_{1}(x):=
\max_{u\in \Real{n},\|u\|_\infty\le 1} u^T(P-I)x+\varepsilon g_1(x)\nn
&=&\|(P-I)x\|_1+\varepsilon g_1(x),
\eeai{almost1}
where
\beai
g_1(x)&=&\varepsilon^{-1}\max_{z\in \Real{n},\,\|z\|_\infty\le \varepsilon,\,|z|_j\le \varepsilon_j} z^Tx\nn
&=&\max_{z\in \Real{n},\,\|z\|_\infty\le 1,\,|z|_j\le {\varepsilon_j/\varepsilon}} z^Tx.
\eeai{conic}
Note that \rf{conic} is a strictly feasible conic optimization problem, when dualizing the constraints we come to
\[
g_1(x)=\min_{u+v=x}\big\{\|u\|_\infty+\sum_j {\varepsilon_j\over \varepsilon}|v_j|\big\}.
\]
Along with \rf{almost1} this implies Proposition \ref{spto1}. $\Box$
\par
The proof of Proposition \ref{spto2} is completely analogous.
\subsection{Proof of Proposition \ref{spto22} }
We start with the following lemma, which is a
modification of Theorem 3.1 in \cite{EGL}.
\begin{lemma} \label{lem1}
For $a\in \Real{n}, b\in \Real{m}, \xi\in \Real{n\times m}$
it holds
\[\max_{\|\xi\|_F\leq\varepsilon}\|a+\xi b\|_2=\|a\|_2+\varepsilon \|b\|_2.\]
\end{lemma}
{\em Proof:} Indeed, $\|a+\xi b\|_2\leq \|a\|_2+\varepsilon
\|b\|_2$, while for $a\neq 0, \;\xi=\xi^*=\frac{\varepsilon
ab^T}{\|a\|_2\|b\|_2}$ equality holds and $\|\xi^*\|_F=\varepsilon$
(for $a=0$ one can take $\xi^*=\frac{\varepsilon cb^T}{\|b\|_2},
\|c\|_2=1, c\in \Real{n}$ arbitrary).\;\;\;$\Box$
\par
We can now bound \bse
\lefteqn{\max_{\xi\in \Xi_F} \|(P+\xi)x-x\|_2 =
\max_{P+\xi\in \mathcal{S},\,||\xi\|_F\le \varepsilon} \|(P+\xi)x-x\|}\\
&\leq&\max_{\|\xi\|_F\leq \varepsilon} \|Px-x+\xi x\|_2=\|Px-x\|_2+\varepsilon \|x\|_2,
\ese
where the last equality is due to Lemma \ref{lem1}. \;\;\;$\Box$

\end{document}